\documentclass[12pt,reqno]{amsart}
\usepackage{enumerate, latexsym, amsmath, amsfonts, amssymb, amsthm, color}
\def\pmod #1{\ ({\rm{mod}}\ #1)}
\def\Z{\Bbb Z}

\def\l{\left}
\def\r{\right}
\def\bg{\bigg}
\def\({\bg(}
\def\){\bg)}
\def\t{\text}
\def\f{\frac}

\def\ls{\leqslant}
\def\gs{\geqslant}

\def\eq{\equiv}

\def\da{\delta}

\def\Proof{\noindent{\it Proof}}

\theoremstyle{plain}
\newtheorem{theorem}{Theorem}

\newtheorem{lemma}{Lemma}
\newtheorem{corollary}{Corollary}

\theoremstyle{definition}

\theoremstyle{remark}

 \vspace{4mm}

\begin{document}

\hbox{Contrib. Discrete Math. 19 (2024), no.\,3, 258--263.}
\medskip

\title
[{A new trigonometric identity with applications}]
{A new trigonometric identity with applications}

\author
[Zhi-Wei Sun and Hao Pan] {Zhi-Wei Sun and Hao Pan}

\address{(Zhi-Wei Sun) School of Mathematics, Nanjing
University, Nanjing 210093, People's Republic of China}
\email{zwsun@nju.edu.cn}

\address {(Hao Pan) School of Applied Mathematics, Nanjing University of Finance and Economics,
 Nanjing 210046, People's Republic of China}
\email {haopan79@zoho.com}

\keywords{Trigonometric identity, roots of unity, Bernoulli polynomials
\newline \indent 2020 {\it Mathematics Subject Classification}. Primary 05A19, 33B10; Secondary 11B68.
\newline \indent The first author is the corresponding author. The two authors were supported by the Natural Science Foundation of China (grants 12371004 and 12071208 respectively).}

\begin{abstract} In this paper we obtain a new curious identity involving trigonometric functions. Namely, for any positive odd integer $n$, we prove that
$$\sum_{k=1}^n(-1)^k(\cot kx)\sin k(n-k)x=\f{1-n}2,$$
which is equivalent to the identity
$$\sum_{k=1}^n(-1)^kU_{n-k}(\cos kx)=-\f{n+1}2,$$
where $U_m(z)$ stands for the $m$th Chebyshev polynomial of the second kind.
As a consequence, for any positive odd integer $n$ and positive integer $m$, we obtain the identity
$$\sum_{k=1}^n(-1)^kk^{2m}B_{2m+1}\left(\frac{n-k}2\right)=0,$$
where $B_j(x)$ denotes the Bernoulli polynomial of degree $j$.
\end{abstract}
\maketitle

\section{Introduction}
\setcounter{theorem}{0}
\setcounter{equation}{0}

Let $\Z^+$ denote the set of all positive integers.
 J.-C. Liu and F. Petrov \cite[(2.11)]{LP} showed that
if $\omega=e^{2\pi i/(3n+2)}$ with $n\in\Z^+$ then
\begin{equation}\label{1.1}\sum_{k=1}^{2n+1}\f{(-1)^k\omega^{k(3k+1)/2}}{1-\omega^{3k}}=-\f{n+1}2,
\end{equation}
which has the equivalent form (cf. \cite[(2.17)]{LP})
\begin{equation}\label{1.2}
\sum_{k=1}^{2n+1}\l(\f{y^{k}}{1+y^{3k}}+\f{(-y)^k}{1-y^{3k}}\r)=-n-1,
\end{equation}
where $y=e^{2\pi i/(6n+4)}$.
Motivated by this, Z.-W. Sun \cite{S} conjectured that
if $m,n\in\{2,3,\ldots\}$ and $\da\in\{0,1\}$, then for any primitive $(m(n-\da)-(-1)^{\da})$-th root
of unity $\zeta$, we have the identity
\begin{equation}\label{1.3}\text{Re}\(\sum_{k=1}^{n-1}\l(\f{\zeta^k}{1+\zeta^{km}}-(-1)^{n+\da}\f{(-\zeta)^k}{1-\zeta^{km}}\r)\)=(-1)^{n-1}\l\lfloor\f n2\r\rfloor.
\end{equation}
This was confirmed by Nemo and Sun in the cases $\delta=0$ and $\delta=1$ respectively; see \cite{S} for the detailed proofs.

Inspired by the above work, we establish the following new result.

\begin{theorem}\label{Th1.1} Let $n$ be any positive odd integer.
Then, for any complex number $q$ with $|q|=1$ and $q^k\not=1$ for all $k=1,\ldots,n$, we have
\begin{equation}\label{1.4}{\rm Re}\(\sum_{k=1}^n\frac{(-1)^kq^{-k(n-k)/2}}{1-q^k}\)=-\f{n+1}4.
\end{equation}
Equivalently, we have the trigonometric identity
\begin{equation}\label{1.5}\sum_{k=1}^{n}(-1)^kU_{n-k}(\cos kx)=-\f{n+1}2,
\end{equation} where $x$ is a real number, $U_m(z)$ is the $m$-th Chebyshev polynomial of the second kind, defined by $U_m(\cos\theta)=(\sin\, (m+1)\theta)/\sin\theta$.
\end{theorem}

\begin{corollary}\label{Cor1.1}
Suppose that $n$ is a positive odd integer and $m$ is a positive integer. Then we have
\begin{equation}\label{1.6}
\sum_{k=1}^n(-1)^kk^{2m}B_{2m+1}\l(\f{n-k}2\r)=0,
\end{equation}
where $B_j(x)$ denotes the Bernoulli polynomial of degree $j$.
\end{corollary}

With the help of Theorem \ref{Th1.1}, we obtain the following result.

\begin{theorem}\label{Th1.2} Let $l,m,n\in\Z^+$ with $l\eq m\pmod2$ and $n\eq1\pmod2$.
Then, for any primitive $(mn+l)$-th root of unity $\zeta$, we have
\begin{equation}\label{1.7}{\rm Re}\(\sum_{k=1}^n\f{\zeta^{k(km+l)/2}}{1-\zeta^{km}}\)=-\f{n+1}4.
\end{equation}
\end{theorem}

Applying Theorem \ref{Th1.2} with $l=1$ and $m=3$, we immediately get the following consequence.

\begin{corollary}\label{Cor1.2} Let $n$ be a nonnegative integer and let $\zeta$ be a primitive $(6n+4)$-th root of unity.
Then
\begin{equation}\label{1.8}\sum_{k=1}^{2n+1}\f{\zeta^{k(3k+1)/2}}{1-\zeta^{3k}}=-\f{n+1}2.
\end{equation}
\end{corollary}

It is interesting to compare our \eqref{1.8} with Liu and Petrov's \eqref{1.1}. Actually, we first found \eqref{1.8} motivated by \eqref{1.1} and then discovered the more general Theorem \ref{Th1.2} and related Theorem \ref{Th1.1}.

We are going to prove Theorem \ref{Th1.1} in the next section, and show Corollary \ref{Cor1.1} and Theorem \ref{Th1.2}
in Section 3.

\section{Proof of Theorem 1.1}
\setcounter{theorem}{0}
\setcounter{equation}{0}

\begin{lemma}\label{Lem2.1} Let $n$ be a positive odd integer, and let $z$ be any complex number. Then
\begin{equation}\label{2.1}\sum_{1\ls k\ls n\atop 0\ls j<(n-k)/2}(-1)^kz^{k(2j+k-n)}=0.
\end{equation}
\end{lemma}
\Proof. Let $\sigma$ denote the left-hand side of \eqref{2.1}. Then, by changing the order of summation, we get
\begin{align*}\sigma=&\ \sum_{j=0}^{(n-3)/2}\sum_{k=1}^{n-2j-1}(-1)^kz^{k(2j+k-n)}
\\=&\ \sum_{j=0}^{(n-3)/2}\sum_{l=1}^{n-2j-1}(-1)^{n-2j-l}z^{(n-2j-l)(2j+(n-2j-l)-n)}
\\=&\ (-1)^n\sum_{j=0}^{(n-3)/2}\sum_{l=1}^{n-2j-1}(-1)^{l}z^{l(2j+l-n)}=-\sigma
\end{align*}
and hence $\sigma=0$. \qed

\medskip
\noindent{\it Proof of Theorem 1.1}. Write $q=e^{2ix}=\cos 2x+i\sin 2x$ with $x$ real, and let
$L$ denote the sum in \eqref{1.4}. Then
\begin{align*}L=&\sum_{k=1}^n(-1)^k\f{\cos k(n-k)x-i\sin k(n-k)x}{1-\cos 2kx-i\sin 2kx}
\\=&\sum_{k=1}^n(-1)^k\f{(1-\cos 2kx+i\sin 2kx)(\cos k(n-k)x-i\sin k(n-k)x)}{(1-\cos 2kx)^2-(i\sin 2kx)^2}
\end{align*}
and hence
\begin{equation*}\begin{aligned}\text{Re}(L)=&\ \sum_{k=1}^n(-1)^k\frac{(1-\cos 2kx)\cos k(n-k)x+(\sin 2kx)\sin k(n-k)x}{2-2\cos 2kx}
\\=&\ \f12\sum_{k=1}^n(-1)^k\cos k(n-k)x+\f12\sum_{k=1}^n(-1)^k(\cot kx)\sin k(n-k)x.
\end{aligned}\end{equation*}
It follows that
\begin{equation}\label{2.2}\begin{aligned}2\,\text{Re}(L)=&\sum_{k=1}^n(-1)^k\f{(\sin kx)\cos k(n-k)x+(\cos kx)\sin k(n-k)x}{\sin kx}\\=&\sum_{k=1}^n(-1)^k\f{\sin k(n+1-k)x}{\sin kx}=\sum_{k=1}^n(-1)^kU_{n-k}(\cos kx).
\end{aligned}
\end{equation}
Thus \eqref{1.4} is equivalent to \eqref{1.5}.

Set $z=e^{ix}$. Then
\begin{align*}2\,\text{Re}(L)=&\ \sum_{k=1}^n(-1)^k\f{z^{k(n+1-k)}-z^{-k(n+1-k)}}{z^k-z^{-k}}
\\=&\ \sum_{k=1}^n(-1)^k\sum_{j=0}^{n-k}(z^k)^j(z^{-k})^{n-k-j}
=\sum_{k=1}^n(-1)^k\sum_{j=0}^{n-k}z^{k(2j+k-n)}.
\end{align*}
For each $k=1,\ldots,n$, clearly
\begin{align*}\sum_{(n-k)/2<j\ls n-k}z^{k(2j+k-n)}=\sum_{0\ls s<(n-k)/2}z^{k(2(n-k-s)+k-n)}=\sum_{0\ls s<(n-k)/2}z^{-k(2s+k-n)}.
\end{align*}
Thus
\begin{align*}2\,\text{Re}(L)=&\ \sum_{k=1\atop 2\mid n-k}^n(-1)^kz^{k(2(n-k)/2+k-n)}
\\&\ +\sum_{k=1}^n(-1)^k\sum_{0\ls j<(n-k)/2}(z^{k(2j+k-n)}+z^{-k(2j+k-n)})
\\=&\ \sum_{r=0}^{(n-1)/2}(-1)^{n-2r}+\sum_{1\ls k\ls n\atop 0\ls j<(n-k)/2}(-1)^k(z^{k(2j+k-n)}+z^{-k(2j+k-n)}).
\end{align*}
Combining this with Lemma \ref{Lem2.1}, we obtain that
$$2\,\t{Re}(L)=(-1)^n\f{n+1}2=-\f{n+1}2$$
and hence \eqref{1.4} follows.

The proof of Theorem \ref{Th1.1} is now complete. \qed

\section{Proofs of Corollary \ref{Cor1.1} and Theorem \ref{Th1.2}}
\setcounter{theorem}{0}
\setcounter{equation}{0}

Recall that the Bernoulli numbers $B_0,B_1,B_2,\ldots$ are given by
$$\f x{e^x-1}=\sum_{k=0}^\infty B_k\f{x^k}{k!}\ \ \ (0<|x|<2\pi).$$

\medskip
\noindent{\it Proof of Corollary \ref{Cor1.1}}. Note that
\begin{align*}&\sum_{k=1}^n(-1)^k\cos k(n-k)x
\\=&\ (-1)^n+\sum_{k=1}^{(n-1)/2}((-1)^k+(-1)^{n-k})\cos k(n-k)x=-1
\end{align*}
since $n$ is odd. Combining this with \eqref{2.2}, we see that \eqref{1.5} has the following equivalent form:
\begin{equation}\label{cot}\sum_{k=1}^n(-1)^k(\cot kx)\sin k(n-k)x=\f{1-n}2.
\end{equation}

It is well known that
$$
\cot x=\sum_{j=0}^\infty\frac{(-1)^j2^{2j}B_{2j}x^{2j-1}}{(2j)!}\ \ \ (0<|x|<\pi)$$
(cf. \cite[p.\,232]{IR}) and
$$\sin x=\sum_{j=0}^\infty\frac{(-1)^jx^{2j+1}}{(2j+1)!}.
$$
So,  by \eqref{cot} we have
\begin{align*}
&\frac{1-n}{2}
=\sum_{k=1}^n(-1)^kx^{2m}\sum_{j=0}^m\frac{(-1)^j2^{2j}B_{2j}k^{2j-1}}{(2j)!}\cdot
\frac{(-1)^{m-j}(k(n-k))^{2m-2j+1}}{(2m-2j+1)!}
\end{align*}
whenever $0<|x|<\pi/n$. Comparing the coefficients of $x^{2m}$ in the both sides of the above equality, we obtain
$$
\sum_{k=1}^n(-1)^kk^{2m}\sum_{j=0}^m\frac{2^{2j}B_{2j}}{(2j)!}\cdot\frac{(n-k)^{2m-2j+1}}{(2m-2j+1)!}=0,
$$
which is equivalent to the desired identity (1.6).
\qed

\medskip
\noindent{\it Proof of Theorem \ref{Th1.2}}. Clearly $L=mn+l$ is even. For $k=1,\ldots,n$, we have
$$\zeta^{k(km+l)/2}=\zeta^{k(L-m(n-k))/2}=(-1)^k\zeta^{-mk(n-k)/2}.$$
Thus
$$\sum_{k=1}^n\f{\zeta^{k(km+l)/2}}{1-\zeta^{km}}=\sum_{k=1}^n\f{(-1)^k(\zeta^m)^{-k(n-k)/2}}{1-(\zeta^m)^k}.$$
Note that
$$L_0:=\f L{\gcd(L,m)}>\f{mn}{\gcd(L,m)}\gs n$$
and $q=\zeta^m$ is a primitive $L_0$-th root of unity. Applying Theorem \ref{Th1.1} we see that the real part of
$$\sum_{k=1}^n\f{\zeta^{k(km+l)/2}}{1-\zeta^{km}}=\sum_{k=1}^n\f{(-1)^kq^{-k(n-k)/2}}{1-q^k}$$
is $-(n+1)/4$. This concludes the proof of Theorem \ref{Th1.2}. \qed

\medskip

\nocite{*}
\bibliographystyle{amsplain}

\end{document}